# About the classification of trigonometric splines

Denysiuk V.P. Dr of Phys-Math. sciences, Professor, Kiev, Ukraine

National Aviation University

kvomden@nau.edu.ua

## Annotation

One of the possible variants of the classification of trigonometric interpolation splines is considered, depending on the chosen convergence factors, the distribution of signs of the basis functions and the interpolation factors. The concept of crosslinking and interpolation grids is introduced; these grids can either match or not match. The proposed classification is illustrated by an example.

## Introduction

The approximation, respectively, of the representation of a known or unknown function through many specific functions can be regarded as a central topic of analysis; such special functions are well defined, easily calculated, and have certain analytical properties [1]. Special functions are often algebraic and trigonometric polynomials, exponential functions [2], polynomials [3] and trigonometric [4] splines, etc. In applications, time is often used as an argument for the functions under research $t$, $t \in [0,T]$; with that in mind, we will denote the functions like $f(t)$.

In many cases, the function is investigated $f(t)$, which is generally at least continuous, has a complex analytical representation, or such representation is unknown at all. In such cases, the function is discretized, that is, replaced by a finite sequence of its instantaneous values at the nodes of some, as a rule, equidistant grid given by $[0,T]$. Then, the task is to approximate this sequence in one way or another by a linear combination of special functions of a particular class. In this approach, it is important to choose a class of special functions that minimizes the deviation of the function $f(t)$ from an approach function built on discrete data.

It is clear that without limiting the class of functions $f(t)$, making such a choice is difficult. We will consider the case when functions $f(t) \in W_v^r$, where $W_v^r$ - is a class of periodic functions having absolutely continuous derivative order $r-1$ ($r = 1, 2, ...$), and the derivative is a function of bounded variation. It is known [5,6] that the best apparatus for approximating class functions $W_v^r$ are simple polynomial splines. The theory of such splines is well-developed (see, e [3], [7], [8], [9]) and etc. The main disadvantage of polynomial splines, in our opinion, is that they have a piecewise structure; this leads to the fact that in practice mainly used third-degree splines that are stitched from pieces of third-degree algebraic polynomials. In addition, such a spline structure greatly limits their application in many computational mathematics problems.

By Schoenberg [10], it was suggested that trigonometric splines are cross-linked from trigonometric polynomials of a certain degree and have the same disadvantage as polynomial splines - a piecewise structure.

In [4], [11], [12], another principle was proposed for constructing trigonometric splines, representing them by uniformly converging trigonometric series (Fourier series), whose coefficients have a certain order of descending. The undeniable advantage of such splines is that they are represented by a single expression throughout the whole interval of the function.

One of the properties of introduced trigonometric splines is that at certain parameter values they coincide with simple polynomial splines [11]; so in this case, all the results of the approximation estimates, are obtained for polynomial splines can be transferred to trigonometric splines. Moreover, the well-developed theory of Trigonometric Fourier series gives reason to expect new results in this field.

An important task of trigonometric spline theory is to classify these splines.

## The purpose of the work.
Development of one of the possible approaches to the classification of trigonometric splines

## The main part.

In [4], interpolation trigonometric splines were considered, in the construction of which were used type of functions

$$C_k^{(I)}(\nu,r,N,t) = \nu_k(r)\cos kt + \sum_{m=1}^{\infty}(-1)^{mI}\left[\nu_{mN+k}(r)\cos(mN+k)t + \nu_{mN-k}(r)\cos(mN-k)t\right];$$

$$S_k^{(I)}(\nu,r,N,t) = \nu_k(r)\sin kt + \sum_{m=1}^{\infty}(-1)^{mI}\left[\nu_{mN+k}(r)\sin(mN+k)t - \nu_{mN-k}(r)\sin(mN-k)t\right], \quad (1)$$

convergence factors $\nu_k(r)$, in descending order $O(k^{-(1+r)})$, and interpolation factors

$$H^{(I)}(r,k) = \nu_k(r) + \sum_{m=1}^{\infty}(-1)^{mI}\left[\nu_{mN+k}(r) + \nu_{mN-k}(r)\right], \quad (2)$$

which provide interpolation of setpoints in the nodes of the grids $\Delta_N^{(I)}$, ($I=0,1$), where $t_j^{(0)} = \frac{2\pi}{N}(j-1)$, and $t_j^{(1)} = \frac{\pi}{N}(2j-1)$, and $N = 2n+1$, ($n = 1,2,...$), and for which their polynomial analogs exist (or can be constructed). However, it is not difficult to build other types of trigonometric interpolation splines that may be interesting in various interpolation problems; the task of classifying these splines is important. It is clear that the task of classifying trigonometric splines is not monosemantic and involves several approaches. One such approach is considered in this paper.

We introduce the function

$$c_k^{(I)}(\nu,r,N,t) = \nu_k(r)\cos kt \pm \sum_{m=1}^{\infty}(-1)^{mI}\left[\nu_{mN+k}(r)\cos(mN+k)t \pm \nu_{mN-k}(r)\cos(mN-k)t\right];$$

$$s_k^{I}(\nu,r,N,t) = \nu_k(r)\sin kt \pm \sum_{m=1}^{\infty}(-1)^{mI}\left[\nu_{mN+k}(r)\sin(mN+k)t \pm \nu_{mN-k}(r)\sin(mN-k)t\right], \quad (3)$$

which we will later call the basic functions of trigonometric splines, convergence factors $\nu_k(r)$, in descending order $O(k^{-(1+r)})$, and interpolation factors

$$hc^{(I)}(r,k) = \nu_k(r) \pm \sum_{m=1}^{\infty}(-1)^{mI}\left[\nu_{mN+k}(r) \pm \nu_{mN-k}(r)\right]$$

$$hs^{(I)}(r,k) = \nu_k(r) \pm \sum_{m=1}^{\infty}(-1)^{mI}\left[\nu_{mN+k}(r) \pm \nu_{mN-k}(r)\right] \quad (4)$$

It is clear that the basis functions (3) and the interpolation factors (4) differ from (1) and (2) by the distribution of the signs before the sums and the signs inside these sums.

As follows from (3), (4), in the general case trigonometric splines depend on the chosen types of convergence factors $\nu_k(r)$, the distribution of signs of the basis functions, the grid of crosslinking and the interpolation grid. Then the crosslink grid index we will mark as $I_1$, and the index of the interpolation grid $I_2$, ($I_1, I_2 = 0,1$). Note that by the term crosslink grid we means a grid that crosslinks polynomial analogs (if any of them exists) of trigonometric splines. It will be shown below that for fixed convergence factors, the crosslink grid determines the basis functions $C_k^{(I_1)}(\nu,r,N,t)$ and $S_k^{(I_1)}(\nu,r,N,t)$; the interpolation grid is determined by factors $hc^{(I_1,I_2)}(r,k)$ і $hs^{(I_1,I_2)}(r,k)$, which are defined in this way

$$hc^{(I_1,I_2)}(r,k) = \nu_k(r) \pm \sum_{m=1}^{\infty}(-1)^{m(I_1-2I_1I_2+I_2)}\left[\nu_{mN+k}(r) \pm \nu_{mN-k}(r)\right],$$

$$hs^{(I_1,I_2)}(r,k) = \nu_k(r) \pm \sum_{m=1}^{\infty}(-1)^{m(I_1-2I_1I_2+I_2)}\left[\nu_{mN+k}(r) \pm \nu_{mN-k}(r)\right]. \quad (5)$$

Note that the crosslink grid and the interpolation grid can both coincide and not coincide.

Our classification of trigonometric splines consists of the following steps..
1. Select the type of convergence factors $\nu_k(r,N)$;

2. Select the type of basis functions that are defined by the distribution of characters of these functions.
3. Determine the type of crosslink grid;
4. Determine the type of interpolation grid;
5. Determine the parameter $r$; note that we will consider several values of this parameter.

Consider this approach for the classification of trigonometric splines more detailed.

We choose the type of convergence factors. For the role of such factors let's choose, for example, factors $\psi(r,k)/k^{1+r}$, where $\psi(r,k)$ - some bounded function of their arguments, and some Furier coefficients of some class functions $W_v^r$ [12], which are in descending order $O(k^{-(1+r)})$.

We choose the type of distribution signs of basic functions. To do this, we draw up a table of distribution signs of these functions

| A1 | A2 | A3 | A4 |
|----|----|----|----|
| B1 | B2 | B3 | B4 |
| C1 | C2 | C3 | C4 |
| D1 | D2 | D3 | D4 |

the elements of which, are tables of signs and are defined as follows:

$$A1=\begin{pmatrix}+ & +\\ + & -\end{pmatrix};\ A2=\begin{pmatrix}+ & +\\ + & +\end{pmatrix};\ A3=\begin{pmatrix}+ & -\\ + & +\end{pmatrix};\ A4=\begin{pmatrix}+ & -\\ + & -\end{pmatrix};$$

$$B1=\begin{pmatrix}- & +\\ + & -\end{pmatrix};\ B2=\begin{pmatrix}- & +\\ + & +\end{pmatrix};\ B3=\begin{pmatrix}- & -\\ + & +\end{pmatrix};\ B4=\begin{pmatrix}- & -\\ + & -\end{pmatrix};$$

$$C1=\begin{pmatrix}- & +\\ - & -\end{pmatrix};\ C2=\begin{pmatrix}- & +\\ - & +\end{pmatrix};\ C3=\begin{pmatrix}- & -\\ - & +\end{pmatrix};\ C4=\begin{pmatrix}- & -\\ - & -\end{pmatrix};$$

$$D1=\begin{pmatrix}+ & +\\ - & -\end{pmatrix};\ D2=\begin{pmatrix}+ & +\\ - & +\end{pmatrix};\ D3=\begin{pmatrix}+ & -\\ - & +\end{pmatrix};\ D4=\begin{pmatrix}+ & -\\ - & -\end{pmatrix}$$

It is clear that each element of the matrix corresponds to two factors $hc^{(I_1,I_2)}(r,k)$ and $hs^{(I_1,I_2)}(r,k)$, which are selected from the conditions of crosslinking and interpolation at the nodes of the grids $\Delta_N^{(I)}$.

According to the proposed classification, trigonometric splines will be referred to as follows
$$St^{(I_1,I_2)}(v,E,r,N,t),$$

Where index $I_1$ defines the crosslink grid, and the index $I_2$ defines an interpolation grid ($I_1, I_2 = 0,1$), parameter determines the type of convergence of selected factors, the parameter $E$ determines the distribution type, parameter $r$ ($r = 1, 2, ...$) determines the smoothness of the spline, parameter $N$ ($N = 2n+1$, $n = 1, 2, ...$) determines the number of nodes of the interpolation uniform grid specified on the interval $[0, 2\pi)$, amd $t$ - is an argument of spline.

To illustrate this classification, let us consider an example. Let on the interval $[0, 2\pi)$ set some grids $\Delta_N^{(I)} = \{t_j^{(I)}\}$, ($t_j^{(0)} = \frac{2\pi}{N}(j-1)$, $t_j^{(1)} = \frac{\pi}{N}(2j-1)$). Also, on this interval, given some periodic continuous function $f(t)$ and let the known values $\{f(t_j^{(I_2)})\}_{j=1}^N = \{f_j^{(I_2)}\}_{j=1}^N$ of this function in the nodes of the grid $\Delta_N^{(I_2)}$. As you know, the coefficients $a_0^{(I_2)}$, $a_k^{(I_2)}$, $b_k^{(I_2)}$, ($k = 1, 2, ..., n$), are calculated by the formulas $a_0^{(I_2)} = \frac{1}{N}\sum_{j=1}^N f_j^{(I_2)}$;

$$a_k^{(I_2)} = \frac{1}{N}\sum_{j=1}^N f_j^{(I_2)} \cos kt_j^{(I_2)};\ \ b_k^{(I)} = \frac{1}{N}\sum_{j=1}^N f_j^{(I)} \sin kt_j^{(I)} \tag{6}$$

Let's put $N = 9$, $\{f_j^{(I_2)}\}_{j=1}^9 = \{3,1,3,2,4,1,3,1,2\}$.

1. We choose the convergence factor.

$$v_k(r,N) = \left[ \frac{\sin\left(\frac{\alpha}{2}k\right)}{k} \right]^{1+r},$$

where $\alpha = 2\pi/N$; it is clear that the order of decreasing of this factor $O(k^{-(1+r)})$. Note that the role of the parameter $\alpha$ requires separate study and is not represented here. 2. Select an element of the character distribution table. We will restict ourselves to only the element A1, for which due to selected convergence factor $v_k(r,N)$ their polynomial analogues exist (or can be constructed). Basic features correspond to this element

$$c_k^{(I_1)}(v,r,N,t) = v_k(r)\cos kt + \sum_{m=1}^{\infty}(-1)^{mI_1}\left[v_{mN+k}(r)\cos(mN+k)t + v_{mN-k}(r)\cos(mN-k)t\right];$$

$$s_k^{(I_1)}(v,r,N,t) = v_k(r)\sin kt + \sum_{m=1}^{\infty}(-1)^{mI_1}\left[v_{mN+k}(r)\sin(mN+k)t - v_{mN-k}(r)\sin(mN-k)t\right],$$

interpolation factors

$$hc^{(I_1,I_2)}(r,k) = v_k(r) + \sum_{m=1}^{\infty}(-1)^{m(I_1-2I_1I_2+I_2)}\left[v_{mN+k}(r) + v_{mN-k}(r)\right],$$

$$hs^{(I_1,I_2)}(r,k) = v_k(r) + \sum_{m=1}^{\infty}(-1)^{m(I_1-2I_1I_2+I_2)}\left[v_{mN+k}(r) - v_{mN-k}(r)\right];$$

and trigonometric spline

$$St^{(I_1,I_2)}(v,A1,r,N,t) = \frac{a_0^{(I_2)}}{2} + \sum_{k=1}^{\frac{N-1}{2}}\left[a_k^{(I_2)}\frac{c_k^{(I_1)}(v,r,N,t)}{hc^{(I_1,I_2)}(r,k)} + b_k^{(I_2)}\frac{s_k^{(I_1)}(v,r,N,t)}{hs^{(I_1,I_2)}(r,k)}\right],$$

which is stitched at the nodes of grid with the index $I_1$ and interpolates the trigonometric polynomial at nodes of grid with the index $I_2$.

3-4. We consider it appropriate to observe both variants of crosslinking grids and interpolation grids.

5. To illustrate, let 's consider several values of a parameter $r$, putting it on $r = 0,1,3$. Note that in the case $r = 0$ the convergence of the trigonometric series that given by the trigonometric splines is not uniform; however, this case is also worth considering.

Here are graphs of this spline for some values of parameter $r$ and different values of grid indexes $I_1$, $I_2$, having regard to the following remarks.

1. Nodes of grid $\Delta_N^{(1)}$ are located between the nodes of the grid $\Delta_N^{(0)}$ and on graphs 1-4 are not displayed.

2. When an element of the character table is selected as an element $A1$, interpolation factors $hc^{(I_1,I_2)}(r,k)$ i $hs^{(I_1,I_2)}(r,k)$ are the same. Taking into account that such a situation does not always take place, we consider it appropriate to give both factors.

3. Charts 1-4 instead of complex spline designations $St^{(I_1,I_2)}(v,A1,r,N,t)$ simplified type designations apply $St1(r,t)$ - $St4(r,t)$, which are correspond to different combinations of index values $I_1, I_2$.

Let's now consider the trigonometric spline $St^{(I_1,I_2)}(v,A1,r,N,t)$.

It is clear that when $I_1 = 0$, $I_2 = 0$ we get;

$$c_k^{(0)}(v,r,N,t) = v_k(r)\cos kt + \sum_{m=1}^{\infty}\left[v_{mN+k}(r)\cos(mN+k)t + v_{mN-k}(r)\cos(mN-k)t\right];$$

$$s_k^{(0)}(v,r,N,t) = v_k(r)\sin kt + \sum_{m=1}^{\infty}\left[v_{mN+k}(r)\sin(mN+k)t - v_{mN-k}(r)\sin(mN-k)t\right],$$

$$hc^{(0,0)}(r,k) = v_k(r) + \sum_{m=1}^{\infty}\left[v_{mN+k}(r) + v_{mN-k}(r)\right]$$

$$hs^{(0,0)}(r,k) = \nu_k(r) + \sum_{m=1}^{\infty}\left[\nu_{mN+k}(r) + \nu_{mN-k}(r)\right]$$

and trigonometric spline $St^{(0,0)}(\nu, A1, r, N, t)$

$$St^{(0,0)}(\nu, A1, r, N, t) = \frac{a_0^{(0)}}{2} + \sum_{k=1}^{\frac{N-1}{2}}\left[a_k^{(0)}\frac{c_k^{(0)}(\nu, r, N, t)}{hc^{(0,0)}(r,k)} + b_k^{(0)}\frac{s_k^{(0)}(\nu, r, N, t)}{hs^{(0,0)}(r,k)}\right].$$

Graphs of this spline for some values of the parameter $r$ is given on pic.1.

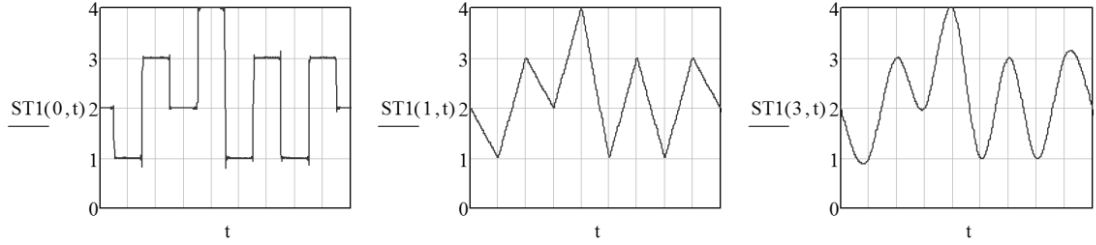

Pic.1 Graphs of spline $St^{(0,0)}(\nu, A1, r, N, t)$ for $r = 0, 1, 3$. Crosslinking and interpolation are carried out in the nodes of the grid $\Delta_N^{(0)}$.

For spline $St^{(0,1)}(\nu, A1, r, N, t)$ while $I_1 = 0$, $I_2 = 1$, we get;

$$c_k^{(0)}(\nu, r, N, t) = \nu_k(r)\cos kt + \sum_{m=1}^{\infty}\left[\nu_{mN+k}(r)\cos(mN+k)t + \nu_{mN-k}(r)\cos(mN-k)t\right];$$

$$s_k^{(0)}(\nu, r, N, t) = \nu_k(r)\sin kt + \sum_{m=1}^{\infty}\left[\nu_{mN+k}(r)\sin(mN+k)t - \nu_{mN-k}(r)\sin(mN-k)t\right],$$

$$hc^{(0,1)}(r,k) = \nu_k(r) + \sum_{m=1}^{\infty}(-1)^m\left[\nu_{mN+k}(r) + \nu_{mN-k}(r)\right],$$

$$hs^{(0,1)}(r,k) = \nu_k(r) + \sum_{m=1}^{\infty}(-1)^m\left[\nu_{mN+k}(r) + \nu_{mN-k}(r)\right],$$

and trigonometric spline $St^{(0,1)}(\nu, A1, r, N, t)$

$$St^{(0,1)}(\nu, A1, r, N, t) = \frac{a_0^{(1)}}{2} + \sum_{k=1}^{\frac{N-1}{2}}\left[a_k^{(1)}\frac{c_k^{(0)}(\nu, r, N, t)}{hc^{(0,1)}(r,k)} + b_k^{(1)}\frac{s_k^{(0)}(\nu, r, N, t)}{hs^{(0,1)}(r,k)}\right].$$

The graphs of this trigonometric spline are shown on pic. 2.

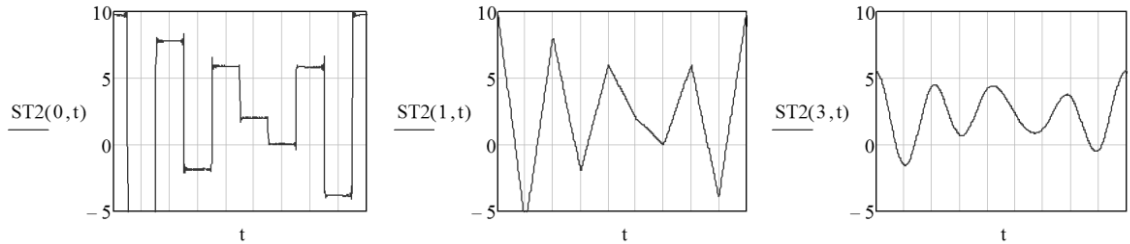

Pic.2 Graphs of spline $St^{(0,1)}(\nu, A1, r, N, t)$ for $r = 0, 1, 3$. Stitching is performed in the nodes of the grid $\Delta_N^{(0)}$, and interpolation at the nodes of the grid $\Delta_N^{(1)}$

For spline $St^{(1,0)}(\nu, A1, r, N, t)$ we get;

$$c_k^{(1)}(\nu, r, N, t) = \nu_k(r)\cos kt + \sum_{m=1}^{\infty}(-1)^m\left[\nu_{mN+k}(r)\cos(mN+k)t + \nu_{mN-k}(r)\cos(mN-k)t\right];$$

$$s_k^{(1)}(\nu,r,N,t) = \nu_k(r)\sin kt + \sum_{m=1}^{\infty}(-1)^m\left[\nu_{mN+k}(r)\sin(mN+k)t - \nu_{mN-k}(r)\sin(mN-k)t\right],$$

$$hc^{(1,0)}(r,k) = \nu_k(r) + \sum_{m=1}^{\infty}(-1)^m\left[\nu_{mN+k}(r) + \nu_{mN-k}(r)\right],$$

$$hs^{(1,0)}(r,k) = \nu_k(r) + \sum_{m=1}^{\infty}(-1)^m\left[\nu_{mN+k}(r) + \nu_{mN-k}(r)\right],$$

and trigonometric spline $St^{(1,0)}(\nu,A1,r,N,t)$

$$St^{(1,0)}(\nu,A1,r,N,t) = \frac{a_0^{(0)}}{2} + \sum_{k=1}^{\frac{N-1}{2}}\left[a_k^{(0)}\frac{c_k^{(1)}(\nu,r,N,t)}{hc^{(1,0)}(r,k)} + b_k^{(0)}\frac{s_k^{(1)}(\nu,r,N,t)}{hs^{(1,0)}(r,k)}\right].$$

The graphs of this trigonometric spline are given in pic. 3.

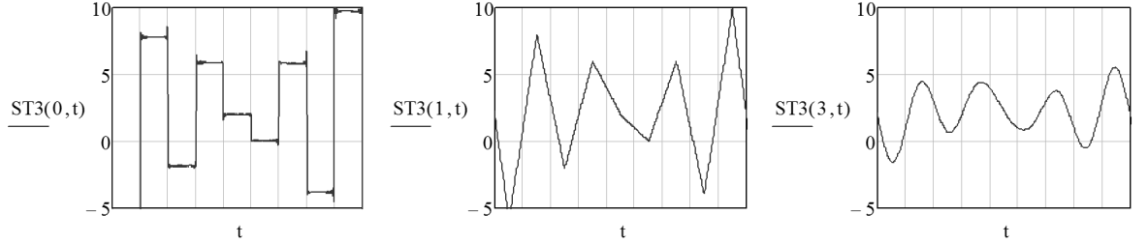

Pic.3 Graphs of spline $St^{(1,0)}(\nu,A1,r,N,t)$ for $r=0,1,3$. Stitching are carried out in the nodes of the grid $\Delta_N^{(1)}$, and interpolation at the nodes of the grid $\Delta_N^{(0)}$ /

For spline $St^{(1,1)}(\nu,A1,r,N,t)$;

$$c_k^{(1)}(\nu,r,N,t) = \nu_k(r)\cos kt + \sum_{m=1}^{\infty}(-1)^m\left[\nu_{mN+k}(r)\cos(mN+k)t + \nu_{mN-k}(r)\cos(mN-k)t\right];$$

$$s_k^{(1)}(\nu,r,N,t) = \nu_k(r)\sin kt + \sum_{m=1}^{\infty}(-1)^m\left[\nu_{mN+k}(r)\sin(mN+k)t - \nu_{mN-k}(r)\sin(mN-k)t\right],$$

$$hc^{(1,1)}(r,k) = \nu_k(r) + \sum_{m=1}^{\infty}\left[\nu_{mN+k}(r) + \nu_{mN-k}(r)\right],$$

$$hs^{(1,1)}(r,k) = \nu_k(r) + \sum_{m=1}^{\infty}\left[\nu_{mN+k}(r) + \nu_{mN-k}(r)\right],$$

and trigonometric spline $St^{(1,1)}(\nu,A1,r,N,t)$

$$St^{(1,1)}(\nu,A1,r,N,t) = \frac{a_0^{(1)}}{2} + \sum_{k=1}^{\frac{N-1}{2}}\left[a_k^{(1)}\frac{c_k^{(1)}(\nu,r,N,t)}{hc^{(1,1)}(r,k)} + b_k^{(1)}\frac{s_k^{(1)}(\nu,r,N,t)}{hs^{(1,1)}(r,k)}\right].$$

The graphs of this trigonometric spline are given in pic.4.

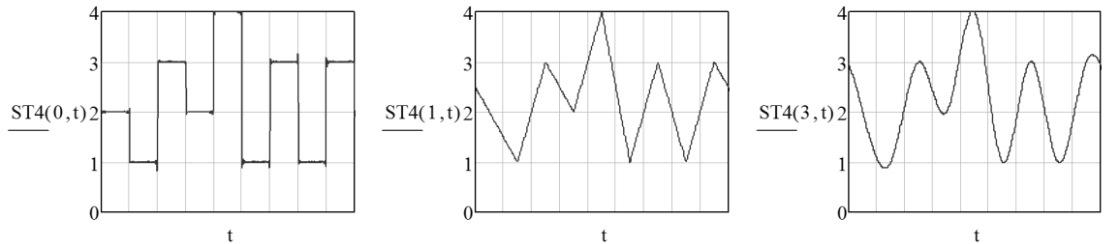

Pic.4 Graphs of spline $St^{(1,1)}(\nu,A1,r,N,t)$ for $r=0,1,3$. Crosslinking and interpolation are carried out in the nodes of the grid $\Delta_N^{(1)}$.

Thus, giving different values to the indexes $I_1$ and $I_2$, you can select crosslinking and interpolation grids based on the needs of the task.

Finally, it should be noted that trigonometric splines are definitely require further investigation.

## Conclusions.

1. One of the variants of classification of trigonometric splines is proposed.

2. Trigonometric splines $St^{(0,0)}(v, A1, r, N, t)$ with odd parameter values $r$ have polynomial analogues - simple polynomial splines of odd degree.

3. Trigonometric splines $St^{(0,1)}(v, A1, r, N, t)$ when the parameter values are even $r$ have polynomial analogues - simple polynomial splines of even degree.

4. Other considered trigonometric splines don't have polynomial analogs ; where appropriate, these analogs can be constructed, at least for cases $r = 0, 1$.

5. Since trigonometric splines are determined by the coefficients of trigonometric interpolation polynomials, well-known FFT (fast Fourier transform) algorithms can be widely used in their construction.

## List of references